 \documentclass[11pt,a4paper]{article}
\usepackage{amsmath}
\usepackage{amsthm}
\usepackage{amssymb}
\usepackage{url}
\usepackage{amsfonts}
\usepackage{setspace}  
\usepackage{anysize}
\usepackage{enumitem}

\usepackage{amscd}        
\usepackage{graphics}
\usepackage{color}
\usepackage[all]{xy}
\usepackage{amsthm}   
\usepackage{mathrsfs}
\usepackage[all]{xy}
\usepackage{textcomp}
\usepackage{fullpage}
\usepackage{mathrsfs}
\usepackage{enumitem}
\usepackage{titlesec}
\usepackage{mathtools}

\titlelabel{\thetitle.\quad}

\titleformat{\subsection}[runin]{\normalfont\bfseries}{\thesubsection.}{3pt}{}
\titleformat{\subsubsection}[runin]{\normalfont\bfseries}{\thesubsubsection.}{3pt}{}

\numberwithin{equation}{section} \theoremstyle{plain}

\theoremstyle{definition}

\theoremstyle{remark}

\def\Z{\mathbb Z}
\def\mcF{\mathcal F}
\def\mcB{\mathcal B}
\def\N{\mathbb N}

\def\Mod{\mbox{-Mod}}

\def\bfC+{\mathbf{C}_+}

\def\C{{\mathcal{C}}}

\def\bfD{\mathbf{D}}
\def\A{\mathcal{A}}
\def\bfC{\mathbf{C}}

\def\A{\mathcal{A}}

\newcommand{\Hom}{\operatorname{Hom}}

\newcommand{\Ext}{\mbox{Ext}}
\newcommand{\id}{\operatorname{id}}

\newcommand{\Proj}{\mbox{Proj}}
\newcommand{\Inj}{\mbox{Inj}}
\newcommand{\Ker}{\mbox{Ker}}
\newcommand{\coker}{\mbox{Coker}}

\newcommand{\GP}{\operatorname{GP}}

\newcommand{\dg}{{\mathit dg\,}}

\newcommand{\Flat}{\mathcal Flat}

\newcommand{\Rep}{\operatorname{Rep}}
\newcommand{\ev}{\operatorname{ev}}
\newcommand{\Cot}{\mathcal Cot}

\setenumerate{label={(\normalfont\arabic*)}} 
\begin{document}

\enlargethispage{\baselineskip}

\title{Completeness of the induced cotorsion pairs in categories of quiver representations}

\date{}

\author{  Sinem Odaba\c{s}{\i}\thanks{e-mail: \texttt{sinem.odabasi@uach.cl}}}

\maketitle
\renewcommand{\theenumi}{\arabic{enumi}}
\renewcommand{\labelenumi}{\emph{(\theenumi)}}

\begin{abstract}
 Given a complete hereditary cotorsion pair $(\A, \mcB)$ in an abelian category $\C$ satisfying  certain conditions,  we study the completeness of the induced cotorsion pairs $(\Phi(\A), \Phi(\A)^{\perp})$  and $(^{\perp}\Psi(\mcB), \Psi(\mcB) )$  in the category $\Rep(Q, \C)$ of $\C$-valued representations of a given quiver $Q$. We show that if $Q$ is left rooted, then the cotorsion pair $(\Phi(\A), \Phi(\A)^{\perp})$ is complete, and if $Q$ is right rooted, then the cotorsion pair $(^{\perp}\Psi(\mcB), \Psi(\mcB) )$ is complete. Besides, we work on the infinite line quiver $A_{\infty}^{\infty}$, which is neither left rooted nor right rooted. We prove that these cotorsion pairs in $\Rep(A_{\infty}^{\infty}, R)$ are complete, as well.

\end{abstract}

\vspace{0.5cm}

\section{Introduction}
Let $(\A, \mcB)$ be a cotorsion pair in an abelian category $\C$, which has  enough $\A$-objects and enough $\mcB$-objects. We now know that there are various ways of lifting such cotorsion pair in $\C$ to a cotorsion pair in, $\bfC(\C)$, the category of chain complexes over $\C$; see for example \cite{Gil08}. Some of these lifts  are:
\begin{enumerate}[label=(\roman*)]
\item $(^{\perp} \bfC(\mcB),\bfC(\mcB) ),$ generated by all disks $D^i(A)$ for every object $A \in \A $.
\item $(\widetilde{\A}, \dg\ \widetilde{\mcB}),$ cogenerated by all spheres $S^i(B)$ for every object $ B \in \mcB$.
\item $(\bfC(\A),\bfC(\A)^{\perp}), $ cogenerated by all disks $D^i(B)$ for every object $B \in \mcB $.
\item $(\dg\ \widetilde{\A}, \widetilde{\mcB}),$ generated by all spheres $S^i(A)$ for every object $ A \in \A$.

\end{enumerate}
The most important application of these induced cotorsion pairs in $\bfC(\C)$ emerges from finding abelian model structures on $\bfC(\C)$, whose homotopy category is the derived category $\bfD(\C)$; see \cite{Hov02} for a detailed treatment about the correspondence between cotorsion pairs and abelian model structures. In \cite{YD15}, the authors proved that whenever  the category $\C$ is (co)complete and  the  cotorsion pair $(\A, \mcB)$ is   complete hereditary,  the so-called  \emph{induced  cotorsion pairs} given in (ii) and (iv) are complete and comptaible, and therefore, yield to such a model structure on $\bfC(\C)$.

Thinking of chain complexes over $\C$ in terms of  $\C$-valued representations of the infinite line quiver with relations, it is natural to ask  how to obtain `quiver  representation version' of these cotorsion pairs, that is, how to  extend a given cotorsion pair in the ground category  $\C$ to a cotorsion pair in the category $\Rep(Q, \C)$ of $\C$-valued representations of any quiver $Q$. Several related results on this problem can be found in the literature for  certain cotorsion pairs in the category $R \Mod$ of (left) $R$-modules, and  for certain quivers, see \cite{EOT04}, \cite{EH99}, \cite{EHHS13}. In \cite{HJ16}, the authors handle the subject in a general framework, working with any abelian category $\C$ and any quiver $Q$, and unify these aforementioned works.

 For every vertex $i $ in a quiver $Q$, under certain conditions, the $i$th  evaluation functor $\ev_i:\ \Rep(Q, \C) \rightarrow \C$  has both the left adjoint $g_i(-)$ and the right adjoint  $f_i(-)$ (see (\ref{f_i})).  Together  with the $i$th stalk functor $s_i:\ \C \rightarrow \Rep(Q, \C)$, which is analogous  to the sphere functor $S^i(-)$, it is proved in \cite[Proposition 7.3]{HJ16} that a cotorsion pair  $(\A, \mcB)$ in an  abelian category $\C$, satisfying certain mild conditions,  leads to the following  cotorsion pairs in $\Rep(Q,\C)$:
\begin{enumerate}[label=(\roman*)]
\item[(i')] $(^{\perp} \Rep(Q,\mcB),\Rep(Q,\mcB) ),$ generated by the class $f_*(\A)$.
\item[(ii')] $(\Phi(\A), \Phi(\A)^{\perp}),$ cogenerated by the class $s_*(\mcB)$.
\item[(iii')] $(\Rep(Q, \A),\Rep(Q, \A)^{\perp} ), $ cogenerated by the class  $g_*(\mcB)$.
\item[(iv')] $(^{\perp} \Psi (\mcB), \Psi(\mcB)),$ generated by the class $s_*(\A)$.
\end{enumerate}
 For the full description  of the aforementioned classes,  we refer to \ref{classes} and Proposition \ref{cotpairrep}. Besides,  they prove that the first (last) two (i')-(ii') ((iii')-(iv')) coincide if the quiver $Q$ is left (right) rooted.

In order to understand better the analogy with the case $\bfC(\C)$, one should be aware  that the $i$th disk functor  $D^i(-)$ is the left adjoint of the $i$th evaluation functor  $\ev_i:\ \bfC(\C) \rightarrow \C$, $\ev_i(X)=X_i$, and besides,  the ($i+1$)th disk functor $D^{i+1}(-)$ is
the right adjoint of $\ev_i$. Furthermore, the intersection of each of the   cotorsion pairs in (i) and (ii) with the subcategory, $\bfC_{-}(\C)$,   of chain complexes  with zeros  in non-negative degrees, are the same, that is, $(^{\perp} \bfC(\mcB) \cap \bfC_{-}(\C),\bfC(\mcB)\cap \bfC_{-}(\C) )= (\widetilde{\A}\cap \bfC_{-}(\C), \dg\ \widetilde{\mcB}\cap \bfC_{-}(\C))$. Dually, the intersection of each of the cotorsion pairs in (iii) and (iv) with the subcategory, $\bfC_{+}(\C)$,   of chain complexes  with zeros  in non-positive degrees, are the same.

From this perspective, we are naturally interested  in answering the question proposed in \cite[Question 7.7]{HJ16} on extending the results given in \cite{YD15} to the cotorsion pairs given in the category $\Rep(Q,\C)$ of $\C$-valued  representations of a quiver $Q$. Namely, in this paper, we study the completeness of the induced cotorsion pairs $(\Phi(\A), \Phi(\A)^{\perp})$ and $(^{\perp} \Psi (\mcB), \Psi(\mcB))$ in $\Rep(Q, \C)$ whenever the cotorsion pair $(\A, \mcB)$ in $\C$ is complete and hereditary.

We now give a summary of the layout   of the paper. In Section 2, we summarize necessary results on quivers and quiver representations. For that, we mostly follow the paper \cite{HJ16}. In Section 3, we provide basic notions and tools  about cotorsion pairs. In Section 4, we prove our main result:

\vskip 0.2in

 \textbf{Theorem \ref{complete}:}
Let $\C$ be an abelian category, and  $(\A,\mcB)$ be a complete hereditary cotorsion pair in $\C$.
\begin{enumerate}
\item If $Q$ is a left rooted quiver,  and if $\C$ has exact $|Q_1^{* \rightarrow i}|$-indexed coproducts for every vertex $i$ in $Q$, then the cotorsion pair $(\Phi(\A), \Phi(\A)^{\perp})$ is complete in $\Rep(Q, \C)$.
\item If $Q$ is a right rooted quiver, and if $\C$ has exact $|Q_1^{i \rightarrow *}|$-indexed products for every vertex $i$ in $Q$, then the cotorsion pair $(^{\perp}\Psi(\mcB), \Psi(\mcB) )$ is complete in $\Rep(Q, \C)$.
\end{enumerate}

\vskip 0.2in

At this point, we note  a  subtle detail in the previous theorem. As clarified in \cite[Remark 4.2]{HJ16}, the cotorsion pairs  $(\Phi(\A), \Phi(\A)^{\perp})$  and $(^{\perp}\Psi(\mcB), \Psi(\mcB) )$ to exist, it isn't necessary to assume the category $\C$ to be  (co)complete and $\C$ to have enough projectives and injectives as stated in \cite[Theorem 7.4]{HJ16}. That's why we just assume $\C$ to have $|Q_1^{* \rightarrow i}|$-indexed coproducts or $|Q_1^{i \rightarrow *}|$-indexed products
for every vertex $i$ in $Q$.

Our approach is inspired from the case of the category $\bfC(\C)$ of chain complexes. Essentially, if the quiver $Q$ is  left rooted, we take  the class  $\Phi(\C)$ as the corresponding one   of acyclic chain complexes. Dually, if $Q$ is right rooted, we consider the class $\Psi(\C)$, instead,  see Lemma \ref{cotpair}.

One should note that if the cotorsion pair $(\A, \mcB)$is generated by a set, and the category $\C$ is a Grothendieck category  with enough $\A$-objects, then the category $\Rep(Q, \C)$ is a Grothendieck category, and the  cotorsion pairs given  in $(i'-iv')$ are also generated by a set, therefore, complete. Even though most of the cotorsion pairs in the theory are known to be generated by a set, there are  examples of complete hereditary cotorsion pairs which are not known to be generated by a set. For instance, if $\C$ is  a (co)complete abelian category with enough projectives, then the trivial cotorsion pair $(\Proj, \C)$ is complete hereditary, but there doesn't exist a  general version of `Kaplansky's Theorem' just as in $R\Mod$, therefore, it is not known to be generated by a set. For another non-trivial example, let $\GP(R)$ denote the class of Gorenstein projective $R$-modules, then  under the conditions given in \cite[Corollary 1]{EIY17}, the cotorsion pair $(\GP(R), \GP(R)^{\perp})$ is complete hereditary but not known to be generated by a set. Having these examples at hand, our result provides a more categorical tool in a (co)complete category  without depending on generation arguments, and extends the class of cotorsion pairs for which the induced cotorsion pairs (ii') and (iv) are complete.

It seems  hard to give an answer in a full generality on the completeness of  the cotorsion pairs $(\Phi(\A), \Phi(\A)^{\perp})$ and $(^{\perp} \Psi (\mcB), \Psi(\mcB))$ for any quiver $Q$. However, in Section 5, we focus on the infinite line quiver $A_{\infty}^{\infty}$, which is neither left rooted nor right rooted, and  we are able prove the following:

\vskip 0.2in

 \textbf{Theorem \ref{infiniteline2}:}
Let $Q$ be the infinite line quiver $A_{\infty}^{\infty}$ and    $(\A,\mcB)$ be  a complete hereditary cotorsion pair in the category $R \Mod$ of  left $R$-modules. Then the cotorsion pairs $(\Phi(\A), \Phi(\A)^{\perp})$ and $(^{\perp} \Psi(\mcB), \Psi(\mcB) )$ in $\Rep(Q, R)$ are complete and hereditary, as well.

\vskip 0.2in

For example, if $\Flat$ and $\Cot$ denote the classes of flat and cotorsion $R$-modules, respectively, then  the pair $(\Flat,\Cot)$
 is known to be a hereditary  cotorsion pair generated by a set, hence, complete. Then, from the construction the cotorsion pair $(^{\perp} \Psi(\Cot), \Psi(\Cot) ))$ in $\Rep(A_{\infty}^{\infty}, R)$ is generated by a set, thus, it is complete. On the other hand, it is not immediate if the cotorsion pair $(\Phi(\Flat), \Phi(\Flat)^{\perp})$ is generated by a set. However, by our result, we can say that it is complete, as well.
\section{Quiver representations}

Throughout this paper, we fix the notation $\C$ for an abelian   category.

\subsection{} Following the notation used in \cite{HJ16}, a quiver, denoted by $Q=(Q_0,Q_1)$, is a directed graph  where $Q_0$ and $Q_1$ are the set of vertices and arrows, respectively. For $i,j \in Q_0$, $Q(i,j)$ denotes the set of paths in $Q$ from $i$ to $j$. For every  $i \in Q_0$, $e_i$ is \emph{ the trivial path}.  If $p \in Q(i,j)$, we write $s(p)=i$ and $t(p)=j$, called its \emph{source} and \emph{target}, respectively.
 Thus, one can think of  $Q$ as a category with the object class $Q_0$ and $\Hom_Q(i,j):= Q(i,j)$.

 For a given vertex $i$  in $Q$, we let
 \begin{equation}
Q_1^{i \rightarrow \ast}:=\{a\in Q_1|\ s(a)=i\}\quad  \textrm{ and }\quad  Q_1^{ \ast \rightarrow i}:=\{a\in Q_1|\ t(a)=i\}.
\end{equation}
From now on, the letter $Q$ will denote a quiver with the sets  of vertices and arrows $Q_0$, $ Q_1$, respectively.

\subsection{} For a given quiver $Q$, the \emph{opposite quiver} $Q^{op}$ is the quiver obtained by reversing arrows of $Q$.

\subsection{} For a given set  $S$, $|S|$ denotes the cardinality of the set $S$. We will say that a category $\C$ has $|Q_1^{* \rightarrow i}|$-\emph{indexed coproducts for every vertex $i$ in a quiver $Q$} if every family $\{C_u\}_{u \in U}$ of objects in $\C$ has the coproduct in $\C$ whenever $|U|  \leq |Q_1^{* \rightarrow i}|$ for some  $i \in Q_0$. Furthermore, if for every  $i \in Q_0$ any $|Q_1^{* \rightarrow i}|$-indexed coproduct of monomorphisms is exact, we say that $\C$ has \emph{exact $|Q_1^{* \rightarrow i}|$-indexed coproducts for every vertex $i $ in $Q$}.

Note that  since in this work we work with abelian categories, if  the  quiver $Q$ is \emph{target-finite}, that is, the set $Q_1^{* \rightarrow i}$ is finite for every $i \in Q_0$, then any abelian category satisfies the aforementioned definitions. So these definitions make sense if  there is a vertex $i$ in the quiver $Q$ such that  $|Q_1^{* \rightarrow i}|$ (  $|Q_1^{i \rightarrow *}| $) is an infinite cardinal.

Dually, we will say that $\C$ has (exact) $|Q_1^{i \rightarrow *}|$-\emph{indexed products for every vertex $i $ in } $ Q$ if $\C^{op}$ has (exact) $|(Q^{op})_1^{* \rightarrow i}|$-\emph{indexed coproducts for every} $i \in Q_0$.

\subsection{Left rooted quiver}\label{leftrootedquiver}\cite{EOT04}, \cite[2.5]{HJ16} For any quiver $Q=(Q_0,Q_1)$, there is a transfinite sequence $\{V_{\alpha}\}_{\alpha}$ of subsets of $Q_0$,  defined as follows:
\begin{itemize}
\item[-] $V_0:=\emptyset$,\quad $V_1:=\{i \in Q_0|\ \textrm{there is no arrow $a$ with } t(a)=i\}$.
\item[-] If $V_\alpha$ is defined for some ordinal $\alpha $,
$$V_{\alpha+1}:=\{ i \in Q_0|\ i \textrm{ is not target of  any arrow a with } s(a) \in Q_0 \backslash \bigcup_{\gamma \leq \alpha} V_{\gamma} \}.$$
\item[-] For a limit ordinal $\beta$, $V_{\beta}:= \bigcup_{\gamma < \beta}V_{\gamma}$.
\end{itemize}

We recall the following results, which are relevant to our purposes and   proved explicitly in \cite{HJ16}.
\subsection{Proposition.}\label{leftrooted-arrow1}\cite[Lemma 2.7]{HJ16} Consider the family $\{V_{\alpha}\}_{\alpha}$ of vertices in a quiver $Q$, defined in \ref{leftrootedquiver}. There is a chain of the form
$$V_1 \subseteq V_2 \subseteq \ldots \subseteq V_{\alpha} \subseteq V_{\alpha+1} \subseteq \ldots \subseteq Q_0.$$

\subsection{Proposition.}\cite[Corollary 2.8]{HJ16}\label{leftrooted-arrow2} If $a : i \rightarrow j $ is an arrow in a quiver $Q$ with  $j \in V_{\alpha+1}$, then $i \in V_{\alpha}$.

\vskip 0.2in

\subsection{Example.} Let $Q$ be the simplest loop    quiver
$ \xymatrix{
.1 \ar@(ul,dl)[]|a }.$
Then $V_{\alpha}= \emptyset$ for every ordinal $\alpha$.

\vskip 0.2in

One can show that if a quiver has an oriented cycle(a nontrivial path with the same source and target), any vertex on it will not belong to  $V_{\alpha}$ for any ordinal $\alpha$. Following the terminology given in \cite{HJ16}, this leads to the so-called \emph{left rooted } quiver, introduced in \cite{EOT04}, in which every vertex belongs to $V_{\alpha}$ for some ordinal $\alpha$: A quiver $Q$ is called \emph{left rooted } if there exists an ordinal $\lambda$ such that $V_{\lambda}=Q_0$. Equivalently, it has no any infinite sequence of arrows $\ldots \bullet \rightarrow \bullet \rightarrow \bullet$, see \cite[Proposition 3.6]{EOT04}.  In a dual manner, a quiver $Q$ is said to be \emph{right rooted} if the  opposite quiver $Q^{op}$ is left rooted.

\subsection{} We let  $\Rep(Q,\C)$ denote  the category of  $\C$-valued functors from $Q$. Equivalently,   $X \in \Rep(Q, \C)$ if and only if $X$ assigns to $i \in Q_0$ an object $X(i) \in \C$, and to every arrow $a: i \rightarrow j$ a morphism $X(a): X(i) \rightarrow X(j)$ in $\C$ with $X(e_i)$ the identity map on $X(i)$. If $p \in Q(i,j)$, then $X(p): X(i) \rightarrow X(j)$ is obtained by composing. An object $X\in \Rep(Q,\C)$ is called  a \emph{$\C$-valued representation} of the quiver $Q$.

A morphism $\lambda: X \rightarrow X'$ between two $\C$-valued representations of $Q$ is just a natural transformation of functors. In fact, $\lambda$ is a family $\{\lambda(i):\ X(i) \rightarrow X' (i)\}_{i \in Q_0}$ of morphisms in $\C$ such that the following diagram commutes
$$
\xymatrix{
X(i) \ar[rr]^{X(a)} \ar[d]_{\lambda (i)} && X(j) \ar[d]^{\lambda(j)}\\
X'(i)\ar[rr]_{X'(a)} && X'(j)}
$$
for every arrow $a: i \rightarrow j$ in $Q$.


\subsection{The evaluation functor.} Given a vertex $i$ in a quiver $Q$, the $i$\emph{th evaluation functor}
$$\ev_i:\ \xymatrix{\Rep(Q,\C) \ar[r] & \C }$$
 assigns to a representation $X \in \Rep(Q,\C)$ the object $X(i)$.

 \subsection{}\label{f_i} Given a quiver $Q$, assume that $\C$ has the coproduct of any family $\{C_u\}_{u \in U}$ where $|U| \leq |Q(i,j)|$ for some $i,j \in Q_0$. Following the notation given in \cite{HJ16}, for a given vertex $i$ in $Q$, we consider the functor $f_i:\ \C \rightarrow \Rep(Q, \C)$ defined as follows: For an object $C \in \C$, the $\C$-valued representation $f_i(C)$ has
$$f_i(C)(j):= \bigoplus_{p \in Q(i,j)} C_p, \quad \textrm{ for all } j \in Q_0 \quad (C_p=C).$$
If $a: j \rightarrow k \in Q_1$ and $p\in Q(i,j)$, then $ap$ is a path from $i$ to $k$. Therefore,  there is a canonical inclusion
$$\iota_{ap}:\xymatrix{C=C_p= C_{ap} \ar@{^{(}->}[r] & f_i(C)(k).} $$
The universal property of coproducts induces  a unique morphism $f_i(C)(a)$, which fits in the following commutative diagram
$$
\xymatrix{C_p\ar@{=}[rr] \ar@{^{(}->}[d]_{\iota_p} && C_{ap} \ar@{^{(}->}[d]^{\iota_{ap}}\\
f_i(C)(j) \ar[rr]_{f_i(C)(a)} && f_i(C)(k).
}
$$

\subsection{}  Given a quiver $Q$, assume that $\C$ has the product of any family $\{C_u\}_{u \in U}$ where $|U| \leq |Q(i,j)|$ for some $i,j \in Q_0$. Then the opposite  category $\C^{op}$ satisfies the condition given in (\ref{f_i}) for the opposite quiver $Q^{op}$. Letting $f_i^{Q^{op}, \C ^{op}}:= f_i$, the functor $g_i$ is defined as $(f_i^{Q^{op}, \C ^{op}})^{op}$.

\subsection{Proposition}\cite[Theorem 3.7]{HJ16} Given a quiver $Q$, suppose that $\C$ has the (co)product of any family $\{C_u\}_{u \in U}$ where $|U| \leq |Q(i,j)|$ for some $i,j \in Q_0$. Let  $i $ be a vertex in  $Q$. The pairs $(f_i,\ev_i)$ and $(\ev_i, g_i)$ are adjoint pairs.

\vskip 0.2in

\subsection{} Note that  all (co)limits in $\Rep(Q, \C)$ are calculated componentwise, hence, the category $\Rep(Q, \C)$ of $\C$-valued representations inherits homological properties of $\C$.

Besides, if $\C$ is a (co)complete abelian category, and if $G$ is a generator (or cogenerator) of $\C$, then the set $\{f_i(G)\}_{i \in Q_0}$ ( $\{g_i(G)\}_{i \in Q_0}$) is a (co)generating set for the category $\Rep(Q,\C)$.

\subsection{The stalk functor.} For a given vertex $i$ in $Q$, the $i$\emph{th stalk functor}
$$s_i:\ \xymatrix{ \C \ar[r] & \Rep(Q,\C)}$$
is defined as
$$s_i(C):=
\left\{
  \begin{array}{ll}
    C, & \hbox{$i=j$;} \\
    0, & \hbox{ $i \neq j$.}
  \end{array}
\right.
$$
with $s_i(a)=0$ for any arrow $a$ in $Q$.

\subsection{} Assume that $\C$ has $|Q_1^{* \rightarrow i}|$-indexed coproducts for every $i \in Q_0$. Let $i$ be a vertex  in $Q$ and  $X \in \Rep(Q,\C)$. By the universal property of coproducts, there is a unique morphism
$$\varphi_i^X:\ \xymatrix{\bigoplus_{a \in Q_1^{* \rightarrow i}} X(s(a)) \ar[r] & X(i)}$$
and $c_i(X):= \coker \ \varphi_i^X$, which  yield a functor
$$c_i:\ \xymatrix{\Rep(Q, \C) \ar[r] & \C }. $$

Dually, if the category $\C$ has $|Q_1^{i \rightarrow *}|$-indexed products for every $i \in Q_0$, the functor $k_i:\ \Rep(Q, \C) \rightarrow \C$ is  defined as $(c_i^{Q^{op}, \C^{op}})^{op}$, that is,
$$k_i(X)= \Ker\ (\xymatrix{X(i) \ar[r]^-{\psi_i^X} & \prod_{a \in Q_1^{i \rightarrow *}} X(t(a))}).$$

\subsection{Proposition.}\cite[Theorem 4.5]{HJ16}
 Let $\C$ be a category with  $|Q_1^{* \rightarrow i}|$-indexed coproducts and $|Q_1^{i \rightarrow *}|$-indexed products for every $i \in Q_0$. Then the pairs $(c_i, s_i)$ and $(s_i, k_i)$ are adjoint pairs.

\vskip 0.2in

\section{Cotorsion Pairs}
 In order to state  our object of interest, we continue recalling  necessary   terminologies and results  on cotorsion pairs.  For a given pair of objects $X,Y$ in a category $\C$ and a non-negative integer $n$, $\Ext^n_{\C}(X,Y)$ denotes the $n$th Yoneda extension class, see \cite[Section III-5]{Mac94}. Even though, it carries an abelian group structure, it may fail to be a set unless $\C$ is an efficient abelian  category, see \cite[Corollary 5.5]{Sto13}. For convenience, we will omit the letter $\C$ in the notation.

For a given class $\A$ of objects of $\C$, we let

$$
\begin{array}{rcl}
  \A ^{\perp} & = & \{X \in \C|\ \Ext^1(A,X)=0 \textrm{ for all } A \in \A\}, \\
  & &\\
  ^{\perp} \A  & = & \{X \in \C|\ \Ext^1(X,A)=0 \textrm{ for all } A \in \A \}.
\end{array}
$$

\subsection{} A pair $(\A, \mcB)$ of classes of objects in $\C$ is said to be a \emph{cotorsion pair} provided that
 $$\A^{\perp}= \mcB \quad  \textrm{ and } ^{\perp} \mcB= \A.$$ The category $\C$ is said to \emph{have enough $\A$-objects} ( $\mcB$-\emph{objects}) if for every object $C \in \C$ there exists an epimorphism $A \twoheadrightarrow C$ (a monomorphism $C \hookrightarrow B$) with $A \in  \A$ ($B \in \mcB$).

Given a class $\mcF$ of objects in $\C$, there are two associated  cotorsion pairs
$$(^{\perp}(\mcF^{\perp}), \mcF^{\perp}), \quad \textrm{ called \emph{cotorsion pair generated by }} \mcF,$$
$$(^{\perp} \mcF, (^{\perp}\mcF)^{\perp}), \quad \textrm{ called \emph{cotorsion pair cogenerated by }} \mcF.$$

\subsection{}Let $(\A, \mcB)$ be a cotorsion pair in $\C$ and  $C$ be an object in $\C$. $C$ is said to \emph{ have  a special $\A$-precover} if there is a short exact sequence

$$ \xymatrix{0 \ar[r] & B \ar[r] & A \ar[r] & C \ar[r] & 0}, \quad \textrm{ with } A\in \A, \textrm{ } B \in \mcB, $$
and to \emph{have a special $\mcB$-preenvelope} if  there is a short exact sequence
$$\xymatrix{0 \ar[r] & C \ar[r] & B' \ar[r] & A' \ar[r] & 0}, \quad \textrm{ with } A' \in \A, \textrm{ } B' \in \mcB.$$

A cotorsion pair $(\A, \mcB)$ is said to \emph{have enough projectives} (\emph{injectives}) if every object of $\C$ has a special $\A$-precover (a special $\mcB$-preenvelope). It is called  \emph{complete} if it has both enough projectives and enough injectives.

The following is a  general form of the so-called  Salce's Lemma, which can be easily proved by using pull-back and push-out arguments.

\subsection{Lemma.}\label{completenesscot}
Let $(\A, \mcB)$ be a cotorsion pair  in a category $\C$.
\begin{enumerate}
\item If $(\A, \mcB)$ has  enough projectives, and if $\C$ has enough $\mcB$-objects, then the cotorsion pair $(\A, \mcB)$ is complete.
\item If $(\A, \mcB)$ has  enough injectives, and if $\C$ has enough $\A$-objects, then the cotorsion pair $(\A, \mcB)$ is complete.
\end{enumerate}

\subsection{} A cotorsion pair $(\A, \mcB)$ in $\C$ is said to be \emph{hereditary} if $\Ext^n(A,B)=0$ for every non-negative integer $n$,   $A \in \A$ and $B \in \mcB$. It is a well-known fact that if the category $\C$ has enough $\A$-objects and $\mcB$-objects,  the following are equivalent:
\begin{enumerate}
\item ($\A, \mcB$) is hereditary.
\item The class $\A$ is closed under kernels of epimorphisms.
\item The class $\mcB$ is closed under cokernels of monomorphisms.
\end{enumerate}

\vskip 0.2in

The following lemma plays a crucial role in our main result. Even though in \cite[Theorem 3.1]{AA02} it is stated for $\C=R \Mod$, its proof is categorical, which permits us to state it  in a more general form for  any abelian category.
\subsection{Lemma.}\label{specialpre}\cite[Theorem 3.1]{AA02} Let $(\A, \mcB)$ be a hereditary cotorsion pair in $\C$ and
$$\xymatrix{0 \ar[r] & C_1 \ar[r] & C \ar[r] & C_2 \ar[r] & 0}$$
\sloppy be an exact sequence in $\C$. If there are short exact sequences $\xymatrix{0 \ar[r] & B_1 \ar[r]  & T_1 \ar[r] & C_1 \ar[r]  & 0 } $ and  $\xymatrix{0 \ar[r] & B_2 \ar[r]  & A_2 \ar[r] & C_2 \ar[r]  & 0 } $ where $B_1, B_2 \in \mcB$ and $A_2 \in \A $, then there is a commutative diagram with exact rows and columns
$$
\xymatrix{&0\ar[d]&0\ar[d]&0\ar[d]&\\
0 \ar[r] & B_1 \ar[r] \ar[d] & B \ar[r] \ar[d] & B_2 \ar[r] \ar[d] & 0 \\
0 \ar[r] & T_1 \ar[r] \ar[d] & T_2 \ar[r] \ar[d] & A_2 \ar[r] \ar[d] & 0 \\
0 \ar[r] & C_1 \ar[r] \ar[d]  & C \ar[r] \ar[d] & C_2 \ar[r] \ar[d] & 0. \\
&0&0&0&
}
$$
where $B \in \mcB$.  If the left column  is a special precover of $C_1$, then  the middle column is  a special $\A$-precover of $C$, as well. The dual statement holds for special $\mcB$-preenvelopes.
\begin{proof}
We just provide a sketch of the proof, which will be needed in the proof of Theorem \ref{infiniteline}. Taking the  pullback of the morphisms $A \twoheadrightarrow C_2$ and $C \twoheadrightarrow C_2$, one obtains a short exact sequence $\mathbb{E}'  \in \Ext^1(A_2, C_1)$. Since the cotorsion pair is hereditary and $B_1 \in \mcB$, $\Ext^1(A_1, T_1) \cong \Ext^1(A_1, C_1 )$, therefore, there exists a unique short exact sequence $\mathbb{E}:0 \rightarrow T_1 \rightarrow T_2 \rightarrow A_2 \rightarrow 0$ whose pushout with the morphism $T_1\twoheadrightarrow C_1$ is $\mathbb{E}'$.

\end{proof}

\subsection{} A \emph{continuous chain} of objects in $\C$ is a directed system  $\{C_{\alpha},\ g_{\alpha \alpha'}\}_{\alpha \leq \alpha' \leq \lambda}$ (indexed by an ordinal $\lambda$) of objects of $\C$ such that for every ordinal $\alpha < \lambda$, the morphism $g_{\alpha, \alpha+1}$ is monic, and for every limit ordinal $\beta \leq \lambda$, $C_{\beta}= \varinjlim_{\alpha < \beta} C_{\alpha}$.

Given a class $\mcF$ of objects of a category $\C$, an object $C \in \C$ is said to \emph{have an $\mcF$-filtration} if there is a continuous chain $\{C_{\alpha}\}_{\alpha \leq \lambda}$ of subobjects of $C$ with $C_0=0$, $C_{\lambda}=C$ and $C_{\alpha+1}/ C_{\alpha} \in \mcF$ for every ordinal $\alpha < \lambda$.

Dually, a \emph{cocontinuous chain} of objects in $\C$ is an inverse system  $\{C_{\alpha},\ g_{\alpha' \alpha}\}_{\alpha \leq \alpha' \leq \lambda}$ (indexed by an ordinal $\lambda$) of objects of $\C$ such that for every ordinal $\alpha < \lambda$, the morphism $g_{\alpha+1, \alpha}$ is epic, and for every limit ordinal $\beta \leq \lambda$, $C_{\beta}= \varprojlim_{\alpha < \beta} C_{\alpha}$.

An object $C \in \C$ is said to \emph{have an $\mcF$-cofiltration} if there is a cocontinuous chain $\{C_{\alpha},\ g_{\alpha' \alpha}\}_{\alpha \leq \lambda}$  of $C$ with $C_0=0$, $C_{\lambda}=C$ and $\Ker\ g_{\alpha+1,\alpha} \in \mcF$ for every ordinal $\alpha < \lambda$.

\vskip 0.2in

The following Lemma, known as Eklof Lemma, plays a crucial role in our main result, and shows that the left class of a cotorsion pair is closed under filtrations while  the right part  is closed under cofiltrations.

\subsection{Eklof Lemma.}\label{Eklof} Let $(\A, \mcB)$ be a cotorsion pair in an abelian category $\C$. Then we have:
\begin{enumerate}
\item If an object $C \in \C$ has an $\A$-filtration, then $C$ belongs to $\A$, as well.
\item If an object $C \in \C$ has a $\mcB$-cofiltration, then $C$ belongs to $\mcB$, as well.
\end{enumerate}
\begin{proof}
Let $\{C_{\alpha}\}_{\alpha \leq \lambda}$ be an $\A$-filtration of $C$. If $\beta \leq \lambda$ is a limit ordinal,  and if  $\mathbb{E}$ is a short exact sequence  in $\Ext^1(C_{\beta}, B)$, where $B \in \mcB$, then consider the short exact sequence $\mathbb{E}g_{\alpha \beta}$  obtained by taking pullback along $C_{\alpha} \rightarrow C_{\beta}$. Apply the arguments given in the proof of  \cite[Theorem 7.3.4]{EJ11}. The arguments for (2) are dual to that of (1). 
\end{proof}

\vskip 0.2in

As claimed, next we  recall the  classes in $\Rep(Q, \C)$ associated to a class  in $ \C$, which are    defined in \cite[Section 7]{HJ16}, and which will be needed  for the construction of cotorsion pairs in $\Rep(Q, \C)$ arising from a cotorsion pair in the ground category $\C$.

 \subsection{}\label{classes} Given a class $\mcF$ of objects in $\C$, we  let
$$
\begin{array}{rcl}
f_*(\mcF) & = & \{ f_i(F) |\ F \in \mcF \textrm{ and } i \in Q_0        \},\\
&&\\
g_*(\mcF) & = & \{ g_i(F) |\ F \in \mcF \textrm{ and } i \in Q_0        \},\\
&&\\
s_*(\mcF) & = & \{ s_i(F) |\ F \in \mcF \textrm{ and } i \in Q_0        \},\\
&&\\
  \Rep (Q, \mcF) & = & \{ X \in \Rep(Q, \C)|\ X(i) \in \mcF \textrm{ for all } i\in Q_0   \}, \\
  &&\\
  \Psi (\mcF) & = & \{ X \in \Rep(Q, \C)|\ \psi_i^X \textrm{ is epic and }  \Ker\ \psi_i^X \in \mcF \textrm{ for all } i\in Q_0   \},\\
  &&\\
 \Phi (\mcF) & = &\{ X \in \Rep(Q, \C)|\ \varphi_i^X \textrm{ is monic and } \coker\ \varphi_i^X \in \mcF \textrm{ for all } i\in Q_0 \}.
\end{array}
$$
as long as they exist.

\subsection{Proposition.}\cite[Proposition 5.2]{HJ16} Let $C \in \C$ and $X \in \Rep(Q, \C)$. Let $i $  be a vertex in a quiver $Q$.
\begin{enumerate}
  \item If the functor $f_i$ exists, for all non-negative integers $n$ there is an isomorphism
$$\Ext^n(f_i(C), X) \cong \Ext^n(C, \ev_i(X)).$$
  \item If the functor $g_i$ exists, for all non-negative integers $n$ there is an isomorphism
$$\Ext^n(X, g_i(\C)) \cong \Ext^n( \ev_i(X), C).$$
\end{enumerate}

The following lemma is a general form of the result given in \cite[Proposition 5.4]{HJ16}.
\subsection{Proposition.}\label{CS} Let $i $ be a vertex in a quiver $Q$. Let $X \in \Rep(Q, \C)$ and $C \in \C$.
\begin{enumerate}
\item Assume that $\C$ has $|Q_1^{* \rightarrow i}|$-indexed coproducts for every $i \in Q_0$. Then there is an inclusion
    $$\Ext^1(c_i(X), C) \hookrightarrow \Ext^1(X, s_i(C)).$$
If $\varphi^X_i$ is a monomorphism, then it is an isomorphism.
\item Assume that $\C$ has $|Q_1^{i \rightarrow *}|$-indexed products for every $i \in Q_0$. Then there is an inclusion
    $$\Ext^1(C, k_i(X)) \hookrightarrow \Ext^1(s_i(C), X).$$
If $\psi^X_i$ is an epimorphism, then it is an isomorphism.
\end{enumerate}

\begin{proof}
The two statements are dual. So we will prove the first statement ($1$). Let $$\mathbb{E}: \xymatrix{0 \ar[r] &C \ar[r]^{\alpha} & T \ar[r]^-{t} & c_i(X) \ar[r] & 0}$$
be a short exact sequence in $\C$. Taking  the pullback of $t$ and $\sigma:\ X(i) \rightarrow c_i(X)$, we get  the following  commutative diagram
$$
\xymatrix{&0\ar[d]&0\ar[d]&\\
&C \ar@{=}[r] \ar[d]_{\alpha'} & C  \ar[d]^{\alpha} & \\
\bigoplus_{a: * \rightarrow i} X(s(a)) \ar[r]^-{\beta} \ar@{=}[d] & T' \ar[r]^{\sigma'} \ar[d]_{t'} & T \ar[r] \ar[d]^t &0 \\
\bigoplus_{a: * \rightarrow i} X(s(a)) \ar[r]_-{\varphi_i^X} & X(i) \ar[d] \ar[r]_{\sigma} & c_i(X) \ar[r] \ar[d] & 0.\\
&0&0&
}
$$
 We let $Y$ be a $\C$-valued representation of $Q$ defined by
$$
Y(j):=\left\{
        \begin{array}{ll}
          T', & \hbox{if } j=i \\
           &  \\
          X(j), & \hbox{if } j \neq i
        \end{array}
      \right.  \qquad  \textrm{ for }j \in Q_0.
$$
For an arrow $a:j \rightarrow k$ in $Q$, we define the morphism  $Y(a):\ Y(j) \rightarrow Y(k)$ as follows:
\begin{itemize}
\item[Case 1:] If $j \neq i$ and $k \neq i$, then $Y(a)=X(a)$.
\item[Case 2:] If $j \neq i$ and $k=i$, then $Y(a):= \iota_a \circ \beta$, where $\iota_a: X(s(a)) \hookrightarrow \bigoplus_{a: * \rightarrow i} X(s(a))$ is the canonical inclusion.
\item[Case 3:] If $j=i$ and $k \neq i$, then $Y(a)= X(a) \circ t'$.
\item[Case 4:] If $j=k=i$, then  $\sigma \circ X(a)=0$, and  by using  the pullback square, there exists a unique morphism $Y(a): T' \rightarrow T'$ which fits in the following commutative diagram

\begin{equation}\label{PB1}
\xymatrix{T'\ar[rd]_{Y(a)}  \ar@/^/[drr]^0 \ar[dd]_{t'}&&\\
& T'\ar[r]_{\sigma'} \ar[dd]^{t'}& T\ar[dd]^t\\
X(i) \ar[rd]_{X(a)} &&\\
& X(i) \ar[r]_{\sigma} & c_i(X).
}
\end{equation}

\end{itemize}
\emph{Claim:} There is a short exact sequence of the form
$$\widetilde{\mathbb{E}}: \xymatrix{0 \ar[r] & s_i(C) \ar[r]^{\widetilde{\alpha}} & Y \ar[r]^{\widetilde{t}} & X \ar[r] & 0}$$
in $\Rep(Q, \C)$ where for a vertex $j$ in $Q$
$$
\widetilde{\alpha}(j):= \left\{
                    \begin{array}{ll}
                      \alpha', & \hbox{ if } j=i \\
                       & \\
                      0 , & \hbox{ if } j \neq i
                    \end{array}
                  \right.
\quad \textrm{ and } \quad
\widetilde{t}(j):=\left\{
                    \begin{array}{ll}
                      t', & \hbox{ if } j=i \\
                       & \\
                      \id , & \hbox{ if } j \neq i.
                    \end{array}
                  \right.
$$
It is easy to verify that $\widetilde{t}:\ Y \rightarrow X$ is a morphism of $\C$-valued representations. As for $\widetilde{\alpha}$, we just need to prove the Case $4$, that is, for an arrow $a:\ i \rightarrow i$, $Y(a) \circ h=0$. From the diagram (\ref{PB1}), we have  $t' \circ Y(a) \circ \alpha'= X(a) \circ t' \circ \alpha'=0$, by uniqueness, $Y(a) \circ \alpha' =0$. Finally, it is obvious that $\widetilde{\mathbb{E}}$ is a short exact sequence in $\Rep(Q, \C)$.

The  assignment to a short exact sequence  $\mathbb{E} \in \Ext^1(c_i(X), C)$  a short exact sequence $\widetilde{\mathbb{E}}$ in $\Ext^1(X, s_i(C))$ is well-defined , therefore, it yields to a homomorphism $\Ext^1(X,s_i(C)) \rightarrow \Ext^1(c_i(X),C)$ of abelian groups. Now we need to show that it is injective. Suppose that the short exact sequence $\widetilde{\mathbb{E}}$ splits and $h=\{h_i\}_{i \in Q_0}$ is a coretraction of $\widetilde{\alpha}$, that is, $h \circ \widetilde{\alpha}= \id$. If $Q_1^{* \rightarrow i}= \emptyset$, then  $\bigoplus_{a: * \rightarrow i} X(s(a))=0$, hence,        $\widetilde{\alpha}(i)=\alpha'=\alpha$ and $h_i \circ \alpha= \id$. If $Q_1^{* \rightarrow i} \neq  \emptyset$, then for every $a \in Q_1^{* \rightarrow i} $, $h_i \circ Y(a)=0$, and therefore $h_i \circ \beta=0$. Since $T= \coker\ \beta$, there exist a unique  morphism $z:\ T \rightarrow C$ such that $z \circ \sigma' = h_i$. Now, $z \circ \alpha=z \circ \sigma' \circ \alpha'=h_i \circ \alpha'= \id$. So the short exact sequence $\mathbb{E}$ splits, as well.

\end{proof}

\vskip 0.2in

The following result shows that under suitable conditions, a cotorsion pair $(\A, \mcB)$ in a category $\C$ gives rises to four canonical cotorsion pairs in $\Rep(Q, \C)$, which appear in the literature when  $\C:=R \Mod$ and $Q$ is   a left or right quiver.

\subsection{Proposition.}\label{cotpairrep}\cite[Theorem 7.4, Theorem 7.9]{HJ16} Let $(\A, \mcB)$ be a cotorsion pair in a (co)complete abelian category $\C$, which  has enough $\A$-objects and $\mcB$-objects. Assume that the cotorsion pair $(\A, \mcB)$ is generated by a class $\A_0$ and cogenerated by a class $\mcB_0$.  Then the following are cotorsion pairs in the category $\Rep(Q, \C)$ of $\C$-valued representations of a quiver $Q$:
\begin{enumerate}
\item $(^{\perp} \Rep(Q, \mcB),\Rep(Q, \mcB))$, generated by $f_*(\A_0)$,
\item $(\Phi (\A), \Phi (\A)^{\perp}),$ cogenerated by $s_*(\mcB_0)$,
\item $( \Rep(Q, \A),\Rep(Q, \A)^{\perp}),$ generated by $g_*(\mcB_0)$,
\item $(^{\perp} \Psi(\mcB),\Psi(\mcB)),$ cogenerated  by $s_*(\A_0)$.
\end{enumerate}
If  $Q$ is left rooted, then the cotorsion pairs given in ($1$-$2$) coincide. If $Q$ is right rooted, then the cotorsion pairs given in ($3$-$4$) coincide.

\subsection{Remark.}\label{remark1} Note that the pairs in Proposition \ref{cotpairrep} are proved to be a cotorsion pair  in \cite[Theorem 7.4]{HJ16} under the condition that $\C$ is a category with enough injectives and projectives, $\Proj  \subseteq \A_0$  and $\Inj \subseteq \mcB_0$.  As a matter of fact, it depends on \cite[Proposition 5.6]{HJ16}, which still holds if $\C$ has enough $\A$-objects and enough $\mcB$-objects.

\subsection{Remark.}\label{remark2} The existence of the cotorsion pairs (2) and (4) in Proposition \ref{cotpairrep} doesn't  depend on the (co)completeness of the category $\C$ so that it can be replaced with the condition that $\C$ has $|Q_1^{* \rightarrow i}|$-indexed  coproducts and $|Q_1^{i \rightarrow *}|$-indexed products for every $i \in Q_0$ just as in Proposition \ref{CS}. Besides, if the quiver $Q$ is left rooted, the transfinite induction done in \cite[Theorem 7.9]{HJ16} would still work, and therefore, we have $\Rep(Q, \mcB) \subseteq \Phi(\A)^{\perp}$. Dually, if $Q$ is right rooted, then $\Rep(Q, \A) \subseteq ^{\perp}\Psi(\mcB)$

\section{Completeness of cotorsion pairs in $\Rep(Q, \C)$}
From now on, we fix the notation  $(\A, \mcB)$ for a cotorsion pair in an abelian category $\C$.  In this section, we aim at proving the completeness of the induced  cotorsion pairs $(\Phi (\A), \Phi (\A)^{\perp})$ and $(^{\perp} \Psi(\mcB),\Psi(\mcB))$ for certain quivers. For that, we need some preparative lemmas.

\subsection{Lemma.}\label{enoughB}
\begin{enumerate}
\item If a cotorsion pair $(\A, \mcB)$ in $\C$ has enough injectives, then the category $\Rep(Q, \C)$ of $\C$-valued representations of a quiver $Q$ has  enough $\Rep(Q, \mcB)$-objects.
\item If a cotorsion pair $(\A, \mcB)$ in $\C$ has enough projectives, then the category $\Rep(Q, \C)$ of $\C$-valued representations of a quiver $Q$ has  enough $\Rep(Q, \A)$-objects.
\end{enumerate}
\begin{proof}
Let $X$ be a $\C$-valued representation of a quiver $Q$. For every $i \in Q_0$, we fix a special $B$-preenvelope
$$\xymatrix{0 \ar[r] & X(i) \ar[r]^{\eta_i} & B_i \ar[r] & A_i \ar[r] & 0}.$$
Since $(\A, \mcB)$ is a cotorsion pair, for every arrow $a:\ i \rightarrow j$ in $Q$, there is a morphism $B(a)$ which fits in the following commutative diagram
$$\xymatrix{X(i)\ar[r]^{\eta_i} \ar[d]_{X(a)} & B_i \ar[d]^{B(a)}\\
X(j) \ar[r]_{\eta_j} & B_j}.$$
We let $B$ be the $\C$-valued representation of $Q$ with $B(i):=B_i$ for every vertex $i$ in $Q$, and  morphisms $B(a)$ for every arrow $a$ in $Q$. By construction, $B \in \Rep(Q, \mcB)$ and the family $\eta:=\{\eta_i\}_{i \in Q_0}$ is a monomorphism from $X$ to $B$.
\end{proof}

\subsubsection{Lemma.}\label{cotpair} 
Let $X$ be a $\C$-valued representation of a quiver $Q$.
\begin{enumerate}
\item If  $Q$ is  a left rooted quiver, and if  $\C$  has enough $\mcB$-objects and $|Q_1^{* \rightarrow i}|$-indexed coproducts for every $i \in Q_0$, then
 there is a short exact sequence
$\xymatrix{0 \ar[r] & B \ar[r] & E \ar[r] & X \ar[r] & 0}$
with $E \in \Phi (\C)$ and $B \in \Rep (Q, \mcB)$.
\item  If $Q$ is  a right rooted quiver, and if   $\C$  has enough $\A$-objects and $|Q_1^{i \rightarrow *}|$-indexed products for every $i \in Q_0$,  then there is a short exact sequence
$\xymatrix{0 \ar[r] & X \ar[r] & E' \ar[r] & A \ar[r] & 0}$
with $E' \in \Psi (\C)$ and $C \in \Rep (Q, \A)$.
\end{enumerate}

\begin{proof}
We will prove the first statement (1), the proof of (2) is dual.

Let $\{V_{\alpha}\}_{\alpha \leq \lambda}$ be the $\lambda$-transfinite sequence of vertices given in (\ref{leftrootedquiver}) with $V_{\lambda}= Q_0$. By transfinite induction, we will construct a continuous inverse $\lambda$-sequence $\{E^{\alpha}\}_{\alpha \leq \lambda}$ in $\Rep(Q, \C)$ which satisfies:
\begin{enumerate}[label=(\roman*)]
\item For every ordinal $0< \alpha < \lambda$ and  $i \in Q_0 \backslash V_{\alpha}$, $E^{\alpha}(i)=X(i)$ and $\Ker (E^{\alpha+1} \twoheadrightarrow E^{\alpha}) \in \Rep(Q, \mcB)$.

\item For every ordinal $\alpha \leq \lambda$ and $i \in V_{\alpha}$,  $\varphi_i^{E^{\alpha}}$ is a monomorphism.

\item  For every $0< \alpha' < \alpha \leq \lambda$, the morphism $E^{\alpha} \twoheadrightarrow  E^{\alpha'}$ is a vertex-wise split epimorphism, and  $E^{\alpha}(i)= E^{\alpha'}(i)$ whenever  $i \notin V_{\alpha} \backslash V_{\alpha'}$.
\end{enumerate}
We set $E^0=0$ and $E^1=X$. Suppose that we have $E^{\alpha}$ for some ordinal $0< \alpha < \lambda$. We now construct $E^{\alpha+1}$ as follows:


For every  $i \in V_{\alpha+1} / V_{\alpha}$, we fix a monomorphism $$\varepsilon_i^{\alpha}:\ \xymatrix{ \bigoplus_{ a \in Q_1^{* \rightarrow i }} E^{\alpha}(s(a)) \ar@{^{(}->}[rr] && B_i^{\alpha+1}},$$
where $B_i^{\alpha+1} \in \mcB$.
Note that if $  Q_1^{* \rightarrow i }= \emptyset$, then both the sum  and $B_i^{\alpha+1}$ are the zero object. If not,  from (\ref{leftrooted-arrow2}), any arrow $a \in Q_1^{* \rightarrow i}$ has the source  $s(a)$  in $ V_{\alpha}$, and therefore, $s(a) \neq i$ and there exists the  canonical inclusion
\begin{equation}\label{sum1}
\varepsilon_i^{\alpha} \circ \iota^{s(a),i}:\ \xymatrix{E^{\alpha}(s(a)) \ar@{^{(}->}[rr]^-{\iota^{s(a),i}} && \bigoplus_{a \in Q_1^{* \rightarrow i}}  E^{\alpha}(s(a)) \ar@{^{(}->}[r]& B_i^{\alpha+1}.}
\end{equation}
By assumption, $i \in Q_0 \backslash V_{\alpha}$, hence, $E^{\alpha}(a):\ E^{\alpha}(s(a)) \rightarrow E^{\alpha}(i)=X(i)$. Using the  universal property of products, there is a canonical morphism
\begin{equation} \label{E1j}
\left(
  \begin{array}{c}
    E^{\alpha}(a) \\
                  \\
    \varepsilon_i^{\alpha} \circ \iota^{s(a),i} \\
  \end{array}
\right):\
\xymatrix{E^{\alpha}(s(a)) \ar[r] & X(i) \oplus  B_i^{\alpha+1}},
\end{equation}
which induces a monomorphism
$$\bigoplus_{a \in Q_1^{* \rightarrow i}} E^{\alpha} (s(a)) \hookrightarrow X(i)\oplus  B_i^{\alpha+1}.$$
 We define $E^{\alpha+1}$ as follows
$$E^{\alpha+1}(i):=\left\{
         \begin{array}{ll}
           E^{\alpha}(i), & i \notin V_{\alpha+1} \backslash V_{\alpha} \\
           &\\
           X(i) \oplus  B_i^{\alpha+1} , &  \hbox{otherwise}.
         \end{array}
       \right.
$$

Let $a :\ j \rightarrow k$ be an arrow in $Q$.
\begin{itemize}
\item[-] If $k \in V_{\alpha+1} \backslash V_{\alpha}$, then $j \in V_{\alpha}$, and  therefore, $E^{\alpha+1}(a)$ is the morphism given in (\ref{E1j}).
\item[-] If $j \in V_{\alpha+1} \backslash V_{\alpha} $, then $k \notin V_{\alpha+1}$, and $E^{\alpha+1}(k)=E^{\alpha}(k)=X(k)$, hence, $E^{\alpha+1}(a):\ E^{\alpha+1}(j) \rightarrow  X(k)$ is the composition of the projection map $E^{\alpha+1}(j) \rightarrow E^ {\alpha}(j)=X(j) $ followed by $E^{\alpha}(a)=X(a)$.
\item[-] For the other cases, $E^{\alpha+1}(a)= E^{\alpha}(a)$.
\end{itemize}
One can easily check that the canonical projection map on each vertex  gives rise to a vertex-wise split epimorphism $E^{\alpha+1} \twoheadrightarrow E^{\alpha}$ with
$$\Ker\ (E^{\alpha+1} \twoheadrightarrow E^{\alpha})= \prod_{i \in V_{\alpha+1} \backslash V_{\alpha}} s_i(   B_i^{\alpha+1}) \in \Rep(Q, \mcB) \subseteq \Phi(\A)^{\perp}.$$
As a  consequence , $E^{\alpha+1}$ satisfies the desired conditions.

If $\beta \leq \lambda $ is a limit ordinal and  $E^{\alpha}$ is constructed for every ordinal $\alpha < \beta$, then we let $E^{\beta}:= \varprojlim_{\alpha < \beta} E^{\alpha}$. In fact, $E^{\beta}$ is of the form
$$E^{\beta}(i):= \left\{
                \begin{array}{ll}
                  E^{\alpha}(i), & i \in V_{\beta} \\
                  &\\
                  X(i), & i \notin V_{\beta}
                \end{array}
              \right.
$$
where $i \in V_{\alpha}$ for some ordinal $\alpha < \beta$ (recall that  $V_{\beta}= \bigcup_{\alpha < \beta} V_{\alpha}$), because for every $i \in V_{\alpha} \subseteq V_{\beta}$, the inverse system $\{E^{\alpha'}(i)\}_{\alpha' < \beta}$ in $\C$ is eventually constant to $E^{\alpha}(i)$. Therefore,  $E^{\beta}$ satisfies the desired conditions.
As a consequence, letting $E:=E^{\lambda}= \varprojlim_{\alpha < \lambda} E^{\alpha}$, we have a short exact sequence
$$\xymatrix{0 \ar[r] & B \ar[r] & E \ar[r] & X \ar[r] & 0 }.$$

To show that $B \in \Rep(Q, \mcB)$, we let $B^{\alpha}:= \Ker ( E^{\alpha}\twoheadrightarrow E^1=X)$. For any ordinals $\alpha \leq \beta \leq \lambda$ there exists a commutative diagram

$$\xymatrix{0 \ar[r] & B^{\beta} \ar[r] \ar[d]_{g_{\alpha \beta}} & E^{\beta} \ar[r] \ar@{->>}[d] & X \ar[r] \ar@{=}[d] & 0\\
0 \ar[r] & B^{\alpha} \ar[r] & E^{\alpha} \ar[r] & X \ar[r] & 0.}$$
Using  pullback arguments, the morphism $g_{\alpha \beta }:\ B^{\beta}  \rightarrow B^{\alpha}$ is a vertex-wise split epimorphism, as well. In other words, $\{g_{\alpha \beta}\}_{\alpha \leq \beta \leq \lambda}$ is a continuous inverse $\lambda$-sequence with $B^{\lambda}=B$. Note that for an ordinal $\alpha < \lambda$
 $$\Ker\ g_{\alpha, \alpha+1} = \Ker\ (E^{\alpha+1} \rightarrow E^{\alpha}) \in \Rep(Q, \mcB).$$
It implies that for every $i \in Q_0$, $B(i)$ has a $\mcB$-cofiltration, by Lemma \ref{Eklof}, $B \in \Rep(Q,\mcB)$.



\end{proof}

\subsubsection{Lemma}\label{monocover2}
Let $E$ and $E'$ be two $\C$-valued representations of a quiver $Q$ which are  in $\Phi(\C)$ and $\Psi(\C)$, respectively. Suppose that   $(\A,\mcB)$ is  a complete hereditary cotorsion pair.
\begin{enumerate}
\item If $Q$ is left rooted and $\C$ has exact $|Q_1^{* \rightarrow i}|$-indexed coproducts for every $i \in Q_0$,  then   $E$ has a special $\Phi(\A)$-precover.
\item   If $Q$ is right rooted and $\C$ has exact $|Q_1^{i \rightarrow *}|$-indexed products for every $i \in Q_0$,  then   $E'$ has a special $\Phi(\mcB)$-preenvelope.

\end{enumerate}
 \begin{proof}
We just prove the first statement because of the duality of arguments.

Let $\{V_{\alpha}\}_{\alpha \leq \lambda}$ be the $\lambda$-transfinite sequence of vertices given in (\ref{leftrootedquiver}) with $V_{\lambda}= Q_0$.  By transfinite induction on ordinals $\alpha \leq \lambda$, we will define  a continuous inverse $\lambda$-sequence of  short exact sequences $\mathbb{E}_{\alpha}: 0 \rightarrow B^{\alpha} \rightarrow A^{\alpha} \rightarrow E^{\alpha} \rightarrow 0$ in $\Rep(Q, \C)$,  which satisfies
 \begin{enumerate}[label=(\roman*)]
 \item $E^{\lambda}=E$
\item If $\alpha < \alpha' \leq \lambda$ and $i \not \in V_{\alpha'} \backslash V_{\alpha}$, then $\mathbb{E}^{\alpha}(i)= \mathbb{E}^{\alpha'}(i)$.
\item for every ordinal $\alpha \leq \lambda$ and $i \in V_{\alpha}$, $\varphi_i^{B^{\alpha}} $ is a monomorphism and $B^{\alpha}(i) \in \mcB$.
\item for every ordinal $\alpha \leq \lambda$ and $i \in V_{\alpha}$, $\varphi_i^{A^{\alpha}} $ is a monomorphism with $ \coker\ \varphi_i^{A^{\alpha}} \in \A$.
\end{enumerate}
Firstly, for every  $i \in Q_0$ we fix a special $\A$-precover $0 \rightarrow B'_i \rightarrow A'_i \rightarrow \coker \ \varphi_i^E \rightarrow 0 $.

 We let $E^{\alpha}$ be the representation
$$ E^{\alpha}(i):=
\left\{
  \begin{array}{ll}
    E(i), & \hbox{$ i \in V_{\alpha}$;} \\
    &\\
    0, & \hbox{\textrm{otherwise}.}
  \end{array}
\right., \quad E^{\alpha}(a):=\left\{
    \begin{array}{ll}
      E(a), & \hbox{$i,j \in V_{\alpha}$;} \\
      &\\
      0, & \hbox{\textrm{otherwise}.}
    \end{array}
  \right.
$$
where  $a: i \rightarrow j\in Q_1$. Clearly $E^0=0$ and $E^{\lambda}=E$. So we let  $\mathbb{E}^0$ be  a short exact sequence of zero representations.  Suppose that we have $\mathbb{E}^{\alpha}$ for some ordinal $ \alpha < \lambda$. We  now construct $\mathbb{E}^{\alpha+1}$ as follows:

 Let $i$ be a vertex  in $ V_{\alpha +1} \backslash V_{\alpha}$. Note that the set $Q_1^{* \rightarrow i}$ may be infinite, so  the coproduct $\bigoplus_{Q_1^{* \rightarrow i}} B^{\alpha}(s(a))$ may
    fail to belong to $\mcB$. Therefore we consider  a special $\mcB$-preenvelope of $\bigoplus_{Q_1^{* \rightarrow i}} B^{\alpha}(s(a))$
$$\xymatrix{0 \ar[r] & \bigoplus_{Q_1^{* \rightarrow i}} B^{\alpha}(s(a)) \ar[r] & B^{\alpha}_i \ar[r] & A_i^{\alpha} \ar[r] &0.
}$$
 Applying  pushout arguments, we have the following commutative diagram
$$
\xymatrix{& 0 \ar[d] & 0 \ar[d] &  &\\
0 \ar[r] & \bigoplus_{Q_1^{* \rightarrow i}} B^{\alpha}(s(a)) \ar[r] \ar[d] &  \bigoplus_{Q_1^{* \rightarrow i}} A^{\alpha}(s(a)) \ar[r] \ar[d] & \bigoplus_{Q_1^{* \rightarrow i}} E^{\alpha}(s(a)) \ar[r] \ar@{=}[d] &0\\
0 \ar[r] &  B^{\alpha}_i \ar[r] \ar[d] & \overline{A}_i^{\alpha} \ar[r] \ar[d] &\bigoplus_{Q_1^{* \rightarrow i}} E^{\alpha}(s(a)) \ar[r]  &0\\
0 \ar[r] & A^{\alpha}_i \ar@{=}[r] \ar[d] & A^{\alpha}_i \ar[d]  & &\\
& 0&0&&
}$$
 Since the class $\A$ is closed under coproducts and extensions, $\overline{A}_i^{\alpha} \in \A$. Note that the first row is exact because  $\C$ has exact $|Q_1^{* \rightarrow i}|$-indexed coproducts for every $i \in Q_0$.    Applying Lemma \ref{specialpre}, we have the following commutative diagram
\begin{equation}\label{1}
\xymatrix{&0\ar[d]&0\ar@{-->}[d]&0\ar[d]& \\
0\ar@{-->}[r] & B_i^{\alpha} \ar[d] \ar@{-->}[r] & B_i \ar@{-->}[d] \ar@{-->}[r] & B'_i \ar[d] \ar@{-->}[r]& 0\\
0\ar@{-->}[r] &\overline{A}_i^{\alpha} \ar[d] \ar@{-->}[r] & A_i \ar@{-->}[d] \ar@{-->}[r] & A'_i \ar[d] \ar@{-->}[r]& 0\\
0\ar[r] & \oplus_{a \in Q_1^{* \rightarrow i}} E^{\alpha}(s(a)) \ar[d] \ar[r] & E(i) \ar@{-->}[d] \ar[r] & \coker\ \varphi_i^E \ar[d] \ar[r]& 0\\
&0&0&0&
}
\end{equation}
in which all three columns are special $\A$-precovers. We set $\mathbb{E}^{\alpha+1}(i):= 0 \rightarrow B_i \rightarrow A_i \rightarrow E(i) \rightarrow 0$ if $i \in V_{\alpha+1} \backslash V_{\alpha}$, and $\mathbb{E}^{\alpha+1}(i)= \mathbb{E}^{\alpha}(i)$ if $i \notin V_{\alpha+1} \backslash V_{\alpha}$.

Let  $a:\ j \rightarrow k \in Q_1$.
\begin{itemize}
\item[-] If  $k \in V_{\alpha+1} \backslash V_{\alpha}$, then $j \in V_{\alpha}$ and  $B^{\alpha+1}(a)$ and $A^{\alpha+1}(a)$ are the canonical morphisms
$$B^{\alpha}(j) \hookrightarrow \bigoplus_{a \in Q_{1}^{* \rightarrow k}} B^{\alpha}(s(a)) \hookrightarrow  B^{\alpha}_j \hookrightarrow B_k,$$
$$A^{\alpha}(j) \hookrightarrow \bigoplus_{a \in Q_{1}^{* \rightarrow i}} A^{\alpha}(s(a)) \hookrightarrow \overline{A}^{\alpha}_k \hookrightarrow A_k.$$
\item[-] If $k \in V_{\alpha}$, then $j \in V_{\alpha}$, and  $B^{\alpha+1}(a)= B^{\alpha}(a)$ and $A^{\alpha+1}(a)= A^{\alpha}(a)$.

\item[-] If $k \in Q_0 \backslash V_{\alpha+1}$, then  $B^{\alpha+1}(k)=B^{\alpha}(k)=0$ and $A^{\alpha+1}(k)=A^{\alpha}(k)=0$. So $B^{\alpha+1}(a)=A^{\alpha+1}(a)=0$.
\end{itemize}

Clearly,  $B^{\alpha+1} \in \Rep(Q, \mcB)$, the morphisms $\varphi_i^{A^{\alpha+1}}$ and $\varphi_i^{B^{\alpha+1}}$ are monomorphisms for every   $i \in V_{\alpha+1}$. Just left  to show that for every $i \in V_{\alpha+1}$, $\coker\  (\varphi_i^{A^{\alpha+1}} ) \in  \A$. Due to the construction, it suffices  to show it for $i \in V_{\alpha+1} \backslash V_{\alpha}$. Applying pullback arguments, we have    the following short exact sequence
$$
\xymatrix{0 \ar[r] & A_i^{\alpha} \ar[r] & \coker\ (\varphi_i^{A^{\alpha+1}}) \ar[r] & A'_i \ar[r] & 0
}.$$
So $\coker\ (\varphi_i^{A^{\alpha+1}}) \in \A$.

If $\beta \leq \lambda$ is a limit ordinal, then $\mathbb{E}^{\beta}:= \varprojlim_{\alpha< \beta}\mathbb{E}^{\alpha}$. It is not hard to verify that  for every vertex $i \in V_{\beta}$, there is an ordinal $\alpha < \beta $ such that  $i \in V_{\alpha}$ and $\mathbb{E}^{\beta}(i):= \mathbb{E}^{\alpha}(i)$, and   if $i \not \in V_{\beta}$, $\mathbb{E}^{\beta}(i)=0$.

By transfinite induction, we have a short exact sequence

$$\mathbb{E}:\ \xymatrix{0 \ar[r] & B \ar[r] & A \ar[r] & E \ar[r] & 0  }$$
where $A:=A^{\lambda} \in \Phi(\A)$ and $B:=B^{\lambda} \in \Rep(Q, \mcB) \subseteq \Phi(\A)^{\perp}$.


 \end{proof}
\subsubsection{Theorem}\label{complete}
Suppose  that  $(\A,\mcB)$ is  a complete hereditary cotorsion pair in $\C$.
\begin{enumerate}
\item If $Q$ is left rooted, and if $\C$ has exact $|Q_1^{* \rightarrow i}|$-coproducts for every vertex $i$ in $Q$, then the cotorsion pair $(\Phi(\A), \Phi(\A)^{\perp})$ is complete in $\Rep(Q, \C)$.
\item If $Q$ is right rooted,  and if $\C$ has exact $|Q_1^{i \rightarrow *}|$-products for every vertex $i$ in $Q$, then the cotorsion pair $(^{\perp} \Psi(\mcB), \Psi(\mcB) )$ is complete in $\Rep(Q, \C)$.
\end{enumerate}
 \begin{proof}
  Using  Lemma \ref{cotpair}, Lemma \ref{monocover2} and  classical pullback-pushout arguments, one may easily see that  the cotorsion pair $(\Phi(\A), \Phi(\A)^{\perp})$ has enough projectives, and the cotorsion pair $(^{\perp} \Psi(\mcB), \Psi(\mcB) )$ has enough injectives. By Lemma \ref{completenesscot}, Lemma \ref{enoughB} and Remark \ref{remark1}, they  are complete, as well.

 \end{proof}

\section{A not left nor right rooted  quiver }
In this section, we are interested in enlarging the class of quivers for which  the induced cotorsion pairs  $(\Phi(\A), \Phi(\A)^{\perp})$  and $(^{\perp} \Psi(\mcB), \Psi(\mcB) )$  are complete whenever the cotorsion pair $(\A, \mcB)$ is so. It seems hard to give a complete answer. However, we are able to prove it   for the infinite line quiver $A_{\infty}^{\infty}$.

\subsection{Lemma.}\label{limcoprod}
Let $\C$ be a (co)complete category. Consider a family $\{F^i:\ J \rightarrow \C\}_{i \in I}$ of functors    over a small category $J$.
\begin{enumerate}
\item If the category $\C$ has  exact coproducts, then the canonical morphism
\begin{equation}\label{lim}
\xymatrix{\bigoplus_{i \in I} \varprojlim_{j \in J} F^i \ar[r] &  \varprojlim_{j \in J} \bigoplus_{i \in I} F^i}
\end{equation}
is monic.

\item If the category $\C$ has  exact products, then the canonical morphism
\begin{equation}\label{lim}
\xymatrix{ \varinjlim_{j \in J} \prod_{i \in I}  F^i \ar[r] &  \prod_{i \in I} \varinjlim_{j \in J} F^i}
\end{equation}
is epic.
\end{enumerate}
\begin{proof}
By duality, we just prove the first statement (1).
Note that for every $i \in I$, the limit $\varprojlim F^i= \Ker\ (\xymatrix{\prod_{j \in J} F^i(j) \ar[r]^-{\eta^i} &  \prod_{\lambda:\ j \rightarrow j'} F ^i(j')} )$, see \cite[Chapter IV, Proposition 8.2]{Ste75}. For every finite subset $S \subset I$, there exists the following commutative diagram with exact rows
$$\xymatrix{0 \ar[r] & \varprojlim \bigoplus_{i \in S} F^i(j) \ar[r] \ar[d] &   \prod_{j \in J} \bigoplus_{i \in S} F^i(j) \ar[r] \ar@{^{(}->}[d] & \prod_{\lambda:\ j \rightarrow j'}\bigoplus_{i \in S} F^i(j') \ar@{^{(}->}[d] \\
0 \ar[r] & \varprojlim \bigoplus_{i \in I} F^i(j)  \ar[r]  &   \prod_{j \in J} \bigoplus_{i \in I} F^i(j) \ar[r]  & \prod_{\lambda:\ j \rightarrow j' }\bigoplus_{i \in I}F^i(j') }$$
whose last two columns are  monic because  products preserve monomorphisms, and thus, the first column is monic, as well. Since the set $S$ is finite, the coproduct in the first row can be taken out of the products and hence, it would be of the form
$$\xymatrix{0 \ar[r] &  \bigoplus_{i \in S} \varprojlim  F^i(j) \ar[r]  &  \bigoplus_{i \in S} \prod_{j \in J} F^i(j) \ar[r]^{\oplus \eta^i}  & \bigoplus_{i \in S} \prod_{j \in J} F^i(j').
}$$
Since coproducts preserve monomorphisms by assumption, the morphism  given in (\ref{lim}) is a monomorphism.
\end{proof}

The following lemma is a kind of generalization of Lemma \ref{cotpair}, but one needs to impose certain  conditions on the category $\C$, instead.

\subsection{Lemma. }\label{monocover3} Let $Q$ be a quiver   without loop.  Let $\C$ be a (co)complete category, and  $(\A, \mcB)$ be a  cotorsion pair in $\C$. Consider  a $\C$-valued representation  $X$ of $Q$.
\begin{enumerate}
\item  If  $\C$  has enough $\mcB$-objects and exact coproducts,
  then there is a short exact sequence
$\xymatrix{0 \ar[r] & B \ar[r] & E \ar[r] & X \ar[r] & 0}$
 of $\C$-valued representations of $Q$ with $E \in \Phi (\C)$ and $B \in \Phi(\A)^{\perp}$.

\item  If  $\C$  has enough $\A$-objects and exact products,
  then there is a short exact sequence
$\xymatrix{0 \ar[r] & X \ar[r] & E' \ar[r] & A \ar[r] & 0}$
 of $\C$-valued representations of $Q$ with $E' \in \Psi (\C)$ and $A \in ^{\perp}\Psi(\mcB)$.
\end{enumerate}
\begin{proof}
We can well order the set $Q_0$. Then there is a unique ordinal $\lambda$ which is order isomorphic to $Q_0$. We rename vertices of $Q$ by ordinals $\alpha < \lambda$.  Notice that the set $\N \times (\lambda\cup \{\lambda\})$ with the  lexicographic order is well-ordered, as well. By transfinite induction on $\N \times (\lambda\cup \{\lambda\}) $, we will construct an inverse system $\{E^{(n, \alpha)}\}_{\substack{n \in \Z\\\alpha \leq \lambda}}$ of $C$-valued representations of the quiver $Q$,  which satisfy:
\begin{enumerate}[label=(\roman*)]
\item For every $(n,\alpha) \leq (n, \alpha')$, the morphism $E^{(n, \alpha')}  \rightarrow E^{(n, \alpha)}$ is a vertex-wise split epimorphism with  $E^{(n, \alpha)}(\beta)= E^{(n, \alpha')}(\beta)$ for every ordinal $\beta \leq \alpha$ or $\alpha' < \beta \leq \lambda $.


\item For every $(n,\alpha)$ with $n > 0$,  $\varphi_{\alpha}^{E^{(n,\alpha)}}$ is a monomorphism.

\item If an ordinal  $\alpha' < \lambda$ is the successor ordinal of $\alpha$, then  the morphism $E^{(n, \alpha')}  \rightarrow E^{(n, \alpha)}$ has the kernel in $\Phi(\A)^{\perp}$ for every $n \in \N$.

\end{enumerate}
For every ordinal $\alpha\leq \lambda$, we let $E^{(0,\alpha)}:=X$.    Now we construct $E^{(1,0)}$ as follows:

Consider a monomorphism $\bigoplus_{a \in Q^{* \rightarrow 0}}X(s(a)) \hookrightarrow B^{(1, 0)}$, where  $B^{(1, 0)} \in \mcB$. We define  $E^{(1,0)}$  by
$$E^{(1,0)}(\alpha):=\left\{
         \begin{array}{ll}
           X(\alpha), & \alpha \neq 0 \\
           &\\
           X(0) \oplus  B^{(1, 0)} , & \alpha=0
         \end{array}
       \right.
$$
for every vertex $\alpha < \lambda$, and with the canonical morphisms between vertices. By assumption, $0 \notin Q_1^{* \rightarrow 0}$, so the morphism $\varphi_0^{E^{(1,0)}}$ is a monomorphism. The morphism $E^{(1,0)} \rightarrow X=E^{(0, \lambda)} $, which is defined by the canonical projection on each vertex,   is  an epimorphism which splits vertex-wise,  and has kernel $s_{0}(B^{(1, 0)} )$. In the same manner, if  $\alpha'= \alpha +1 < \lambda$, the $\C$-valued representation $E^{(1, \alpha')}$ of $Q$ is constructed through  $E^{(1, \alpha)}$ and the vertex $\alpha'$.


Now let $\beta < \lambda$ be a limit ordinal. Suppose that $E^{(1, \alpha)}$ is constructed for every ordinal $\alpha < \beta$. Then we let $T^{(1, \beta)}=\varprojlim_{\alpha < \beta}E^{(1, \alpha)}$. Due to the construction, one can easily check that the $T^{(1, \beta )}$  is of the form
\begin{equation}\label{limit1}
T^{(1,\beta)}(\alpha'):=\left\{
         \begin{array}{ll}
           E^{(1, \alpha')}(\alpha') , & \alpha' < \beta \\
           &\\
           X(\alpha') , & \beta \leq \alpha' < \lambda.
         \end{array}
       \right.
\end{equation}
because for every vertex $\alpha' < \lambda$ the sequence $\{E^{(1, \alpha)}(\alpha')\}_{\alpha < \beta}$ is eventually stable. Hence,  there is a splitting epimorphism $T^{(1, \beta)} \rightarrow X=E^{(0, \beta)}$.  In the same way as done above, we proceed the argument in order to   define $E^{(1, \beta)}$ through $T^{(1, \beta )}$ and the vertex $\beta$. Finally we define $E^1:=E^{(1, \lambda)}:= \varprojlim_{\alpha < \lambda}E^{(1, \alpha)}$, which is similar to the form given in (\ref{limit1}), and the canonical vertex-wise split epimorphism $E^1 \rightarrow X$. The family $\{E^{(1, \alpha)}\}_{\alpha \leq \lambda}$ with the morphisms defined forms a cocontinuous $\lambda$-sequence.

For every $\alpha \leq \lambda$,  we let $B^{(1, \alpha) }:= \Ker\ (E^{(1, \alpha)} \rightarrow X)$. Then for every $\alpha \leq \beta \leq \lambda$, there is a commutative diagram with exact rows
$$\xymatrix{0 \ar[r] & B^{(1, \beta)} \ar[r] \ar[d]_{g_{(1,\alpha \beta)}} & E^{(1,\beta)} \ar[r] \ar@{->>}[d] & X \ar[r] \ar@{=}[d] & 0\\
0 \ar[r] & B^{(1,\alpha)} \ar[r] & E^{(1,\alpha)} \ar[r] & X \ar[r] & 0}$$

By pullback arguments, the morphism $g_{(1,\alpha \beta) }:\ B^{(1,\beta)}  \rightarrow B^{(1,\alpha)}$ is a vertex-wise split epimorphism, as well. Therefore, $\{g_{(1,\alpha \beta)}\}_{\alpha \leq \beta \leq \lambda}$ is a cocontinuous inverse $\lambda$-sequence with $\varprojlim_{\alpha < \lambda} B^{(1, \alpha)}=B^{(1,\lambda)}=B^1$. Note that for a successor ordinal $\alpha' =\alpha+1 < \lambda$
 $$\Ker\ g_{(1,\alpha, \alpha')} = \Ker\ (E^{(1,\alpha')} \rightarrow E^{(1,\alpha)})= s_{\alpha'}(B^{(1, \alpha')}) \in \Phi(\A)^{\perp}.$$
In other words,   $\{g_{(1,\alpha \beta)}\}_{\alpha \leq \beta \leq \lambda}$ is a $\Phi(\A)^{\perp}$-cofiltration of $B^{1}$. As a consequence, $B^1 \in \Phi(\A)^{\perp}$.

 For a given integer $n>0 $, using $\C$-valued representation $E^{n-1}:=E^{(n-1, \lambda)}$, we repeat the previous  argument to obtain the family $\{E^{(n, \alpha)}\}_{\alpha \leq \lambda}$, and finally    a vertex-wise split epimorphism $E^{n}:=E^{(n, \lambda)} \rightarrow E^{n-1} $ with the kernel in $\Phi (\A)^{\perp}$.   Finally, we have a cocontinuous inverse system
$$\ldots \rightarrow E^2 \rightarrow E^1 \rightarrow E^0=X$$
whose morphisms between consecutive objects  are vertex-wise splitting epimorphisms and $\Ker\ (E^{n+1} \rightarrow  E^n)=B^{n+1} \in \Phi (\A)^{\perp}$. We denote $E:= \varprojlim_{n \in  \N } E^n$.  By \cite[Lemma 66]{Murfet}, there is a short exact sequence
$$\xymatrix{0 \ar[r] & B \ar[r] & E \ar[r] & X \ar[r] & 0}$$
where $B := \varprojlim_{n \in \Z} B^n$. As done before, one can easily prove that $B$ has  a $\Phi(\A)^{\perp}$-cofiltration $\{B^n\}_{n \in \N}$, therefore $B \in \Phi(\A)^{\perp}$.

Now  we need to prove that $E \in \Phi(\C)$. Let $\alpha$ be a  vertex in $Q$, that is, an ordinal $\alpha < \lambda$. Since  the set $\{(n,\alpha)\}_{n >0}$ is  cofinal  in $\N  \times (\lambda \cup \{\lambda\})$,  $E(\alpha)= \varprojlim_{n >0 }E^{(n, \alpha)}(\alpha)$. Since  $\varphi_{\alpha}^{E^{(n, \alpha)}}$ is monic for every $n >0$, and inverse limits preserves monomorphisms,     the morphism
$$\varprojlim_{n  >0} \varphi_{\alpha}^{E^{(n,\alpha)}}:\ \varprojlim_{n >0} \bigoplus_{a \in Q_1^{*\rightarrow \alpha}} E^{(n,\alpha)}(s(a)) \rightarrow E(\alpha)$$
is a monomorphism, as well. By using the universal property of limits,  the morphism $\varphi_{\alpha}^E$ can be factorized as
$$ \bigoplus_{a \in Q_1^{*\rightarrow \alpha}} E(s(\alpha))=\bigoplus_{a \in Q_1^{*\rightarrow \alpha}}\varprojlim_{n >0} E^{(n,\alpha)}(s(a)) \rightarrow \varprojlim_{n >0} \bigoplus_{a \in Q_1^{*\rightarrow \alpha}} E^{(n,\alpha)}(s(a)) \hookrightarrow E(\alpha)$$
with the canonical morphisms. By Lemma \ref{limcoprod}, the first morphism is monic, as well.

\end{proof}

\subsection{Lemma.}\label{infiniteline}
Let $Q$ be the infinite line quiver
$$A_{\infty}^{\infty}:\ \quad \xymatrix{\ldots\ar[r] & \bullet_{1} \ar[r] & \bullet_0 \ar[r] & \bullet_{-1} \ar[r] & \bullet_{-2} \ar[r]& \ldots, }$$
and $(\A, \mcB)$ be a complete hereditary cotorsion pair in $R \Mod$.
\begin{enumerate}
\item For any $E\in \Phi (R \Mod) \subseteq \Rep(Q, R)$, there is a short exact sequence
$$\xymatrix{0 \ar[r] & B \ar[r] & A \ar[r] & E \ar[r] & 0}$$
with $A \in \Phi (\A)$ and $B \in \Phi (A)^{\perp}$.

\item For any $E'\in \Psi (R \Mod) \subseteq \Rep(Q, R)$, there is a short exact sequence
$$\xymatrix{0 \ar[r] & E' \ar[r] & B' \ar[r] & A' \ar[r] & 0}$$
with $B' \in \Psi (\mcB)$ and $A' \in ^{\perp}\Psi (A)$.
\end{enumerate}
\begin{proof}
We just prove (1), the second one is dual. Consider an object $E \in \Phi (\C)$
$$E:\ \quad  \xymatrix{\ldots\ar[r] & E_{1} \ar@{^{(}->}[r] & E_0 \ar@{^{(}->}[r] & E_{-1} \ar@{^{(}->}[r] & E_{-2} \ar@{^{(}->}[r]& \ldots }.  $$
For every vertex $i$ in $Q$, we fix  a special $\A$-precover $\mathbb{E}_i:\ 0 \rightarrow B_i \rightarrow   A_i \rightarrow   c_i(X) \rightarrow  0.$
Taking the pullback of the morphisms $E_0 \twoheadrightarrow c_0(E)$ and $A_0 \twoheadrightarrow c_0(E)$, and using Lemma \ref{specialpre}, we have the following short exact sequence $0 \rightarrow B^0 \rightarrow A^0 \rightarrow E \rightarrow 0$ of $R$-module valued representations of $Q$

$$
\scalebox{0.8}{
\xymatrix{
B^0:\ \ldots \ar[r] & 0\ar[dd]\ar[r]&0 \ar[dd]\ar[r]&B^0_0 \ar@{^{(}->}[rr]  \ar@{=}[rd]  \ar@{^{(}->}[dd] &&B^0_{-1} \ar@{^{(}->}[rr] \ar@{->>}[rd] \ar@{^{(}->}[dd] && B^0_{-2} \ar@{^{(}->}[dd] \ar@{->>}[rd]\ar@{^{(}->}[r] & \ldots \\
&&& & B_0  \ar@{^{(}->}[dd] && B_{-1}  \ar@{^{(}->}[dd]  & & B_{-2} \ar@{^{(}->}[dd] \\
A^0:\ \ldots \ar[r]& E_2\ar@{=}[dd] \ar@{^{(}->}[r] &E_1 \ar@{^{(}->}[r] \ar@{=}[dd]& T^0_0\ar@{^{(}->}[rr] \ar@{->>}[rd] \ar@{->>}[dd] &&  T^0_{-1} \ar@{^{(}->}[rr] \ar@{->>}[rd] \ar@{->>}[dd]&& T^0_{-2} \ar@{->>}[dd]\ar@{->>}[rd]\ar@{^{(}->}[r] & \ldots \\
&&& & A_0 \ar@{->>}[dd]  && A_{-1}  \ar@{->>}[dd] & & A_{-2}\ar@{->>}[dd] \\
E:\ \ldots \ar[r]&E_2\ar@{^{(}->}[r] &E_1 \ar@{^{(}->}[r] & E_0\ar@{^{(}->}[rr]\ar@{->>}[rd] &&  E_{-1} \ar@{^{(}->}[rr] \ar@{->>}[rd] && E_{-2} \ar@{->>}[rd] \ar@{^{(}->}[r]& \ldots\\
&&& & c_0(E) && c_{-1}(E)  && c_{-2}(E) \\
}}
$$
From the construction just as in  the proof of Lemma \ref{specialpre} , for every integer $j \geq 0 $, the induced short exact sequence of cokernels $ 0 \rightarrow c_j(B^0) \rightarrow c_j(A^0) \rightarrow c_j(E) \rightarrow 0$ is the same as $\mathbb{E}_j$.

In the same manner, for any integer $i \geq 0$ one can obtain a short exact sequence $0\rightarrow B^i \rightarrow A^i \rightarrow E \rightarrow 0$ in $\Rep(Q, R)$
with $B^i, A^i \in \Phi (R \Mod)$,  and the induced short exact sequence $0 \rightarrow c_j(B^i) \rightarrow C_j(A^i) \rightarrow c_j(E) \rightarrow 0 $ is the same as $\mathbb{E}_j$, for every $j \leq i$, and $B^i(j)=0$, $j > i$. One can easily show that $B^i \in \Phi(\mcB)^{\perp}$ because it has a $s_*(\mcB)$-cofiltration.

For every $i \geq 0$, the morphism  $A^{i+1} \twoheadrightarrow E$ has a factorization
$$\xymatrix{A^{i+1} \ar@{->>}[r] \ar@{->>}[rd] & A^i\ar@{->>}[d]\\
& E
}
$$
 such that $\Ker\ (A^{i+1} \twoheadrightarrow A^i )=f_{i+1}(B_{i+1})\in \Phi(\mcB)^{\perp}$.  We will show it for $i=0$. If $j \geq 2$, $A^1(j)=A^0(j)=E_j$. For $j=1$, from the  construction, there is an epimorphism $A^1(1)=T_1^1 \twoheadrightarrow E_1=A^0(0)$ with the kernel $B_1$. Now suppose that for some $j\geq 0$ the epimorphism  $A^1(j)= T^1_j \twoheadrightarrow E_j$ is factorized over the epimorphism $A^0(j)=T^0_j \twoheadrightarrow E_j$, that is, $T^1_j \twoheadrightarrow T^0_j \twoheadrightarrow E_j$ with $\Ker\ ( T^1_j \twoheadrightarrow T^0_j)=B_1$. Since the cotorsion pair $(\A, \mcB)$ is hereditary,  we have a commutative diagram
$$
\scalebox{0.8}{
\xymatrix{
\Ext^1(A_{j-1}, T^1_j) \ar[rr]^-{\cong}\ar@{=}[dd] && \Ext^1(A_{j-1}, T_j^0) \ar[dd]^{\cong}\\
&&\\
\Ext^1(A_{j-1}, T^1_j) \ar[rr]_-{\cong} && \Ext^1(A_{j-1},E_{j}).
}}
$$
Taking the pullback of the morphisms $A_{j-1}\twoheadrightarrow c_{j-1}(E)$ and $E_{j-1}\twoheadrightarrow c_{j-1}(E)$, we obtain  a short exact sequence $\mathbb{F} \in \Ext^1(A_{j-1},E_j)$. Using arguments in  the proof of Lemma \ref{specialpre} and from the previous  diagram, the short exact sequence $0 \rightarrow T_j^0 \rightarrow T_{j-1}^0 \rightarrow A_{j-1} \rightarrow 0$ is the image of the short exact sequence $0 \rightarrow T_j^1 \rightarrow T_{j-1}^1 \rightarrow A_{j-1} \rightarrow 0$, that is, pushout along the morphism $T^1_j \twoheadrightarrow T^0_j$. As a result, there is a morphism  $T^1_{j-1} \twoheadrightarrow T^0_{j-1}$ with the kernel $B_1$. It proves our claim.

Finally , we have an inverse system of short exact sequences of $R$-module representations of the quiver $Q$
$$\xymatrix{
0 \ar[r]& B^{i+1} \ar[r] \ar[d]& A^{i+1} \ar[r] \ar[d] & E \ar[r] \ar@{=}[d] & 0\\
 0 \ar[r]& B^{i} \ar[r] & A^i \ar[r]  & E \ar[r]  & 0.
}$$
Taking the limit of this system, finally we have a short exact sequence $0 \rightarrow B \rightarrow A \rightarrow E \rightarrow 0$. Since $\Ker\ (B^{i+1} \rightarrow B^i)= f_{i+1}(B_{i+1})\in \Phi(\mcB)^{\perp}$, by Lemma \ref{Eklof}$, B \in \Phi(\mcB)^{\perp}$. Besides $A \in \Phi(R \Mod)$ because  limits preserve monomorphisms. Note that $c_j(A^i)=A_j$ for every $j\leq i$,  therefore, $c_j(A)=A_j$ for every integer $j$. As a consequence $A \in \Phi (\A)$.
\end{proof}

\subsection{Theorem.}\label{infiniteline2}
Under the setup as in Lemma \ref{infiniteline}, the cotorsion pairs $(\Phi(\A), \Phi(\A)^{\perp})$ and $(^{\perp}\Psi(\mcB), \Psi(\mcB))$ are complete.

\begin{proof}
By Lemma \ref{monocover3},  Lemma \ref{infiniteline} and pullback-pushout arguments, the cotorsion pairs $(\Phi(\A), \Phi(\A)^{\perp})$ and $(^{\perp}\Psi(\mcB), \Psi(\mcB))$ have enough projectives and injectives, respectively. Note also the category $\Rep(Q, R)$ has enough $\Phi(\A)$-objects and $\Psi(\mcB)$-objects. In fact, for every $X \in \Rep(Q, R)$ we have
$$\xymatrix{X \ar@{^{(}->}[r] &\prod_{i \in Q_0}f_i(B_i) }, \quad \textrm{ and } \quad
\xymatrix{\bigoplus_{i \in Q_0}g_i(A_i) \ar@{->>}[r] & X,}$$
where for every $i \in Q$, $X(i) \hookrightarrow B_i$ is a monomorphism with $B_i \in \mcB$, and $A_i \twoheadrightarrow X(i)$ is an epimorphism with $A_i \in \A$. Then Lemma \ref{completenesscot} can be applied.
\end{proof}


\begin{thebibliography}{99}
\bibliographystyle{alpha}

\bibitem[AA02]{AA02} Akinci, K. \& Alizade, R.  (2002). Special precovers in cotorsion theories. Proceedings of the Edinburgh Mathematical Society, 45(2), 411-420.
\bibitem[EH99]{EH99} Enochs, E. \& Herzog, I. (1999). A homotopy of quiver morphisms with applications to representations.  Canad. J. Math. 51, no:2, 294-308.

\bibitem[EJ11]{EJ11} Enochs, E. E. \& Jenda, O. M. G. Relative Homological Algebra
Volume 1, De Gruyter Expositions in Mathematics 30, Walter de Gruyter, 2000

\bibitem[EOT04]{EOT04} Enochs, E., Oyonarte, L. \& Torrecillas, B. (2004) Flat covers and flat representations of quivers. Comm. Algebra 32, no:4, 1319-1338.
\bibitem[EHHS13]{EHHS13} Eshraghi, H., Hafezi, R., Hosseini, E. \& Salarian, S. (2013). Cotorsion Theory in the category of quiver representations. Journal of Algebra and Its Applications. Vol. 12, No. 6.
 \bibitem[EIY17]{EIY17} Estrada, S., Iacob, A. \& Yeomans, K. (2017). Mediterr. J. Math. 14: 33. https://doi.org/10.1007/s00009-016-0822-5.


\bibitem[Gil08]{Gil08} Gillespie, J. (2008). Cotorsion pairs and degreewise homological model structures. Homology, Homotopy and Applications, vol. 10(1), 283-304.



\bibitem[HJ16]{HJ16}Holm, H. \& Jorgensen, P. (2016).  Cotorsion pairs in categories of quiver representations. Kyoto J. Math. to appear.

\bibitem[Hov02]{Hov02} Hovey, M. (2002). Cotorsion pairs, model category structures, and representation theory. Math. Z. 241, 553–-592.

\bibitem[Mac94]{Mac94} Maclane, S.  Homology. Springer-Verlag Berlin Heidelberg, fourth printing.
\bibitem[Mur06]{Murfet} Murfet, D. (2006). Derived Categories Part I. \url{http://therisingsea.org/notes/DerivedCategories.pdf}.
\bibitem[Ste75]{Ste75} Stenstr\"{o}m, B. (1975). Rings of quotients. Die Grundlehren der mathematischen
Wissenschaften in Einzeldarstellungen, Band 217, Springer-Verlag, New
York, Heidelberg, Berlin.

\bibitem[Sto13]{Sto13}  Stovicek, J. (2013). Exact model categories, approximation theory, and cohomology of quasi-coherent sheaves. Advances in Representation Theory of Algebra. 	arXiv:\ 1301.5206.

\bibitem[YD15]{YD15} Yang, X. \& Ding, N. (2015). On a question of Gillespie. Forum Math. 27, 3205-3231.



\end{thebibliography}
\end{document}